\documentclass[12pt,leqno]{amsart}
\newtheorem*{theorem}{Theorem}
\newtheorem*{corollary}{Corollary}

\begin{document}

\title{A Note on Trimedial Quasigroups}
\author{Michael~K.~Kinyon}
\address{Department of Mathematical Sciences \\
Indiana University South Bend \\
South Bend, IN 46634 USA}
\email{mkinyon@iusb.edu}
\author{J.~D.~Phillips}
\address{Department of Mathematics \& Computer Science \\
Wabash College \\
Crawfordsville, IN 47933 USA}
\email{phillipj@wabash.edu}

\begin{abstract}
The purpose of this brief note is to sharpen a result of Kepka
\cite{kepka1} \cite{kepka2} about the axiomization of the variety
of trimedial quasigroups.
\end{abstract}

\maketitle
\thispagestyle{empty}

A groupoid is \emph{medial} if it satisfies the identity $wx\cdot yz = wy\cdot xz$.
A groupoid is \emph{trimedial} if every subgroupoid generated by $3$
elements is medial. Medial groupoids and quasigroups have also been called
abelian, entropic, and other names, while trimedial quasigroups have also
been called triabelian, terentropic, etc. (See \cite{beneteau}, especially
p. 120, for further background.)

In \cite{kepka1} \cite{kepka2}, Kepka showed that a quasigroup satisfying the following
three identities must be trimedial.
\begin{eqnarray}
xx\cdot yz &=& xy\cdot xz \\
yz\cdot xx &=& yx\cdot zx \\
(x\cdot xx)\cdot uv &=& xu \cdot (xx\cdot v)
\end{eqnarray}

\noindent The converse is trivial, and so these three identities
characterize trimedial quasigroups. Here, we show that, in fact,
(2) and (3) are sufficient to characterize this variety (as a subvariety
of the variety of quasigroups). Note that in the theorem we only assume
left cancellation, not the full strength of the quasigroup axioms.

\begin{theorem}
A groupoid with left cancellation which satisfies (2) and (3) must also satisfy (1).
\end{theorem}

\begin{proof}
$(x\cdot xz)(xx\cdot yz) = (x\cdot xx)(xz\cdot yz) = (x\cdot xx)(xy\cdot zz)
= (x\cdot xy)(xx\cdot zz) = (x\cdot xy)(xz\cdot xz) = (x\cdot xz)(xy\cdot xz)$.
Now cancel.
\end{proof}

In \cite{kepka1} \cite{kepka2}, Kepka showed that the following single identity
characterizes
trimedial quasigroups:
\[
[(xx\cdot yz)]\{ [xy\cdot uu] [(w\cdot ww)\cdot zv]\}
=
[(xy\cdot xz)]\{ [xu\cdot yu] [wz\cdot (ww\cdot v)]\}.
\]
Using the theorem we can sharpen this.

\begin{corollary}
The following identity characterizes trimedial quasigroups:
\[
[(xy\cdot uu)] [(w\cdot ww)\cdot zv] = [(xu\cdot yu)][wz\cdot (ww\cdot v)].
\]
\end{corollary}

\begin{proof}
To obtain (2) set $z = ww$ and use right cancellation. To obtain (3)
set $y = u$ and use left cancellation.
\end{proof}

\smallskip
\noindent\textbf{Acknowledgement.}~Our investigations were aided
by the automated reasoning tool OTTER, developed by McCune \cite{otter}.
We thank Ken Kunen for introducing us to OTTER.

\end{document}